\documentclass[a4paper, 12pt]{article}%
\usepackage{amsmath}

\usepackage{amsfonts}%
\usepackage{amssymb}

\begin{document}

\begin{center}
{\Huge A mean ergodic theorem for actions of amenable quantum groups}

\bigskip

Rocco Duvenhage

\bigskip

\textit{Department of Mathematics and Applied Mathematics}

\textit{University of Pretoria, 0002 Pretoria, South Africa}

\bigskip

2007-9-10
\end{center}

\bigskip

\noindent\textbf{Abstract}

We prove a weak form of the mean ergodic theorem for actions of amenable
locally compact quantum groups in the von Neumann algebra setting.

\bigskip


\section{Introduction}

The following mean ergodic theorem is well-known: Let $G$ be a locally compact
group with right Haar measure $\mu$, and assume that it contains a F\o lner
net $\left(  \Lambda_{\lambda}\right)  $, i.e. a net of Borel sets in $G$ such
that $0<\mu(\Lambda_{\lambda})<\infty$ and $\lim_{\lambda}\mu\left(
\Lambda_{\lambda}\Delta\left(  \Lambda_{\lambda}g\right)  \right)  /\mu\left(
\Lambda_{\lambda}\right)  =0$ for all $g\in G$. Furthermore, let $U_{g}$ be a
contraction on a Hilbert space $H$ such that $U_{g}U_{h}=U_{gh}$ for all
$g,h\in G$, and $G\ni g\mapsto\left\langle U_{g}x,y\right\rangle $ is Borel
measurable for all $x,y\in H$. Take $P$ to be the projection of $H$ onto
$V:=\left\{  x\in\mathfrak{H}:U_{g}x=x\text{ for all }g\in G\right\}  $. Then
\begin{equation}
\lim_{\lambda}\frac{1}{\mu\left(  \Lambda_{\lambda}\right)  }\int
_{\Lambda_{\lambda}}U_{g}xd\mu(g)=Px\tag{1.1}%
\end{equation}
for all $x\in H$. A standard proof for the case $G=\mathbb{Z}$ can be found
for example in \cite{K} and \cite{P}, but it can be extended to the more
general case without much effort.

In this paper we prove a version of this theorem for the action of an amenable
locally compact quantum group on a von Neumann algebra. We use the von Neumann
algebra setting for quantum groups, as developed by Kusterman and Vaes
\cite{KV4} building on earlier work on Kac algebras (see for example \cite{EN}).

In this setting a \textit{locally compact quantum group} is defined to be a
von Neumann algebra $M$ with a unital normal $\ast$-homomorphism
$\Delta:M\rightarrow M\otimes M$ (where $M\otimes N$ denotes the von Neumann
algebraic tensor product of two von Neumann algebras), such that
$(\Delta\otimes\iota_{M})\circ\Delta=(\iota_{M}\otimes\Delta)\circ\Delta$
(where $\iota_{M}$ denotes the identity map on $M)$, and on which there exists
normal semi-finite faithful (n.s.f.) weights $\varphi$ and $\psi$ such that
$\varphi\left(  (\theta\otimes\iota_{M})\circ\Delta(a)\right)  =\varphi
(a)\theta(1)$ for all $a\in\mathcal{M}_{\varphi}^{+}$ and $\psi\left(
(\iota_{M}\otimes\theta)\circ\Delta(a)\right)  =\psi(a)\theta(1)$ for all
$a\in\mathcal{M}_{\psi}^{+}$, for all $\theta\in M_{\ast}^{+}$, where
$M_{\ast}^{+}$ is the positive normal linear functionals on $M$, and
$\mathcal{M}_{\varphi}^{+}=\left\{  a\in M^{+}:\varphi(a)<\infty\right\}  $.
This quantum group is denoted as $\left(  M,\Delta\right)  $. We refer the
reader to the papers \cite{KV1, KV2, KV3} for background and motivation for
this definition. If furthermore there exists a net $\left(  \varphi_{\lambda
}\right)  $ of normal states on $M$ such that $\left\|  \theta\ast
\varphi_{\lambda}-\varphi_{\lambda}\right\|  $ converges to $0$ for all
$\theta\in M_{\ast}$ with $\theta(1)=1$, then we call $(M,\Delta)$
\textit{amenable}; see for example \cite{DQV}. Here $\mu\ast\nu:=(\mu
\otimes\nu)\circ\Delta$ for any $\mu,\nu\in M_{\ast}^{+}$.

An \textit{action} of $(M,\Delta)$ on another von Neumann algebra $A$ is
defined to be a normal injective unital $\ast$-homomorphism $\alpha
:A\rightarrow M\otimes A$ such that $\left(  \iota_{M}\otimes\alpha\right)
\circ\alpha=(\Delta\otimes\iota_{A})\circ\alpha$; see \cite{V}.

Given such an action, we will assume the presence of a normal state $\omega$
on $A$ which is \textit{invariant} under the action, by which we mean that
$(\theta\otimes\omega)\circ\alpha=\omega$ for all normal states $\theta$ on
$M$. In Section 2 we show how to set up the analogue of the integral in (1.1)
for a quantum group action, and in Section 3 we state and prove a mean ergodic
theorem for such actions, however only in a weak form analogous to%
\[
\lim_{\lambda}\left\langle x,\frac{1}{\mu\left(  \Lambda_{\lambda}\right)
}\int_{\Lambda_{\lambda}}U_{g}yd\mu(g)\right\rangle =\left\langle
x,Py\right\rangle
\]
for all $x,y\in H$. Our approach is to set the problem up in a suitable
Hilbert space framework, closely related to that of (1.1), and then to follow
the basic structure of (1.1)'s proof.

We will not need the full force of the theory of locally compact quantum
groups as developed in \cite{KV3, KV4}, and therefore it will be convenient to
formulate our results in an abstract setting incorporating only the concepts
from locally compact quantum groups that we need, modelled on the definitions
discussed above. We will focus on this abstract setting, rather than on
concrete examples.

\section{A suitable integration theory}

In this section we develop the tools and notation that we need in order to
formulate and prove the mean ergodic theorem in the next section. Throughout
this section and the next we will use the following notation: $R$ will be an
arbitrary von Neumann algebra, its unit denoted by $1_{R}$, and its normal
states by $\left(  R_{\ast}^{+}\right)  _{1}$. By $\omega$ we will mean an
arbitrary normal state on a von Neumann algebra $A$. We will denote the GNS
construction of $(A,\omega)$ by $(H,\gamma)$, by which we mean that $H$ is a
Hilbert space and $\gamma:A\rightarrow H$ a linear mapping such that
$\left\langle \gamma(a),\gamma(b)\right\rangle =\omega(a^{\ast}b)$ and with
$\gamma(A)$ dense in $H$.

We remind the reader that we will use the notation $R\otimes A$ to indicate
the von Neumann algebraic tensor product, often written as $R\overline
{\otimes}A$ in the literature. The algebraic tensor product will be written as
$R\odot A$. We will constantly use tensor products of mappings on von Neumann
algebras, and a useful reference regarding this topic is \cite{St}. For
example, if $\theta$ is a normal state on $R$ while $\iota_{A}$ is the
identity map $A\rightarrow A$, then we can define $\theta\otimes\iota
_{A}:R\otimes A\rightarrow A$ as the tensor product of conditional
expectations, in which case $\theta\otimes\iota_{A}$ itself is a conditional
expectation, which is also normal, i.e. $\sigma$-weakly continuous; see
\cite[Section 9]{St}.

We are going to view $R$ as a noncommutative measurable space, and roughly
speaking we will be integrating $A$ valued ``functions'' over $R$.

Note that the integral in (1.1) is an integral of a bounded function
$f:G\rightarrow H$ which can be defined via the Riesz representation theorem
by $\left\langle \int_{\Lambda}fd\mu,x\right\rangle =\int_{\Lambda
}\left\langle f(g),x\right\rangle d\mu(g)$. We now mimic this construction for
$A$ valued ``functions'' on $R$, in other words for elements of $R\otimes A$.

\bigskip

\noindent\textbf{Proposition 2.1. }\textit{Let }$\mu$\textit{ be a normal
positive linear functional on }$R$\textit{. Then there is a unique function }%
\[
\tilde{\mu}:R\otimes A\rightarrow H
\]
\textit{such that }%
\begin{equation}
\left\langle \gamma(d),\tilde{\mu}(T)\right\rangle =\mu\otimes\omega\left(
\left[  1_{R}\otimes d\right]  ^{\ast}T\right)  \tag{2.1}%
\end{equation}
\textit{ for all }$T\in R\otimes A$\textit{ and }$d\in A$\textit{.
Furthermore, }$\tilde{\mu}$\textit{ is linear, }$\left\|  \tilde{\mu}\right\|
\leq\left\|  \mu\right\|  $\textit{ and }%
\[
\left\langle \gamma(d),\tilde{\mu}(T)\right\rangle =\omega\left(  d^{\ast
}\left(  \mu\otimes\iota_{A}\right)  (T)\right)
\]
\textit{for all }$T\in R\otimes A$\textit{ and }$d\in A$\textit{.}

\bigskip

\noindent\textbf{Proof.} For any $T\in R\otimes A$, define the linear
functional%
\[
f_{T}:\gamma(A)\rightarrow\mathbb{C}:\gamma(d)\mapsto\overline{\mu
\otimes\omega\left(  \left[  1_{R}\otimes d\right]  ^{\ast}T\right)  }%
\]
which is indeed well defined, since if $\gamma(d)=0$ we have $\mu\otimes
\omega\left(  \left[  1_{R}\otimes d\right]  ^{\ast}T\right)  =0$ as follows:
First consider any $T=\sum_{j=1}^{n}r_{j}\otimes a_{j}\in R\odot A$, then
\[
\left|  \mu\otimes\omega\left(  \left[  1_{R}\otimes d\right]  ^{\ast
}T\right)  \right|  \leq\sum_{j=1}^{n}\left|  \mu\left(  r_{j}\right)
\right|  \left|  \omega\left(  d^{\ast}a_{j}\right)  \right|
\]
but $\left|  \omega\left(  d^{\ast}a_{j}\right)  \right|  \leq\sqrt
{\omega\left(  d^{\ast}d\right)  }\sqrt{\omega\left(  a_{j}^{\ast}%
a_{j}\right)  }=0$, since $\omega\left(  d^{\ast}d\right)  =\left\|
\gamma(d)\right\|  ^{2}$, therefore $\mu\otimes\omega\left(  \left[
1_{R}\otimes d\right]  ^{\ast}T\right)  =0$. For a general $T\in R\otimes A$
there is a net $T_{\lambda}\in R\odot A$ converging $\sigma$-weakly to $T$,
according to von Neumann's density theorem (see for example \cite[Section
2.4.2]{BR}). Hence $\left[  1_{R}\otimes d\right]  ^{\ast}T_{\lambda}$
converges $\sigma$-weakly to $\left[  1_{R}\otimes d\right]  ^{\ast}T$, but
$\mu\otimes\omega$ is $\sigma$-weakly continuous (i.e. normal), so $\mu
\otimes\omega\left(  \left[  1_{R}\otimes d\right]  ^{\ast}T\right)  =0$.

Clearly $f_{T}$ is linear, and $\left\|  f_{T}\right\|  \leq\left\|
\mu\right\|  \left\|  T\right\|  $ since%
\begin{align*}
\left|  f_{T}\left(  \gamma(d)\right)  \right|   &  \leq\sqrt{\mu\otimes
\omega\left(  \left[  1_{R}\otimes d\right]  ^{\ast}\left[  1_{R}\otimes
d\right]  \right)  }\sqrt{\mu\otimes\omega\left(  T^{\ast}T\right)  }\\
&  \leq\sqrt{\mu\left(  1_{R}\right)  }\sqrt{\omega(d^{\ast}d)}\sqrt{\left\|
\mu\otimes\omega\right\|  \left\|  T^{\ast}T\right\|  }\\
&  =\sqrt{\left\|  \mu\right\|  }\left\|  \gamma(d)\right\|  \sqrt{\left\|
\mu\right\|  }\left\|  T\right\|
\end{align*}
Therefore $f_{T}$ can be linearly extended uniquely to $H$ without changing
its norm. By the Riesz represention theorem and since $\gamma(A)$ is dense in
$H$, there is a unique element $\tilde{\mu}(T)$ in $H$ such that $f_{T}%
(\gamma(d))=\left\langle \tilde{\mu}(T),\gamma(d)\right\rangle $ for all $d\in
A$. Furthermore $\left\|  \tilde{\mu}(T)\right\|  =\left\|  f_{T}\right\|  $.
Hence we obtain a unique function $\tilde{\mu}:R\otimes A\rightarrow H$ such
that (2.1) holds. Clearly $\tilde{\mu}$ is linear and $\left\|  \tilde{\mu
}(T)\right\|  \leq\left\|  \mu\right\|  \left\|  T\right\|  $.

Lastly, for $r\in R$ and $a\in A$ we have
\[
\mu\otimes\omega\left(  \left[  1_{R}\otimes d\right]  ^{\ast}(r\otimes
a)\right)  =\omega(d^{\ast}\left(  \mu\otimes\iota_{A}\right)  (r\otimes a))
\]
hence $\mu\otimes\omega\left(  \left[  1_{R}\otimes d\right]  ^{\ast}T\right)
=\omega(d^{\ast}\left(  \mu\otimes\iota_{A}\right)  T)$ for all $T\in R\odot
A$ by linearity. But again by $\sigma$-denseness, and by $\sigma$-weak
continuity, this extends to all $T\in R\otimes A$. $\square$

\bigskip

We now take this a step further by finding an analogue of the linear operator
$H\rightarrow H:x\mapsto\int_{\Lambda}U_{g}xd\mu(g)$ that appears in (1.1).

\bigskip

\noindent\textbf{Proposition 2.2.}\textit{ Consider the situation in
Proposition 2.1 and furthermore assume that we have a }$\ast$%
\textit{-homomorphism }$\alpha:A\rightarrow R\otimes A$\textit{ which leaves
}$\omega$\textit{ invariant in the following sense:}%
\begin{equation}
(\mu\otimes\omega)\circ\alpha=\mu\left(  1_{R}\right)  \omega\tag{2.2}%
\end{equation}
\textit{for the given }$\mu$\textit{. Then there exists a unique linear
operator }$\tilde{\mu}^{\alpha}:H\rightarrow H$\textit{ such that}%
\[
\tilde{\mu}^{\alpha}\left(  \gamma(a)\right)  =\tilde{\mu}\circ\alpha(a)
\]
\textit{for all }$a\in A$\textit{. Furthermore, }$\left\|  \tilde{\mu}%
^{\alpha}\right\|  \leq\left\|  \mu\right\|  $\textit{, and if }$\alpha
$\textit{ is unital, then }$\left\|  \tilde{\mu}^{\alpha}\right\|  =\left\|
\mu\right\|  $\textit{.}

\bigskip

\noindent\textbf{Proof.} The operator $\tilde{\mu}^{\alpha}$ is well defined
on $\gamma(A)$ since $\tilde{\mu}\circ\alpha(a)=0$ when $\gamma(a)=0$ as we
now show: For any $d\in A$ we have from Proposition 2.1 that%
\begin{align*}
\left|  \left\langle \gamma(d),\tilde{\mu}\circ\alpha(a)\right\rangle \right|
^{2}  &  =\left|  \mu\otimes\omega\left(  \left[  1_{R}\otimes d\right]
^{2}\alpha(a)\right)  \right|  ^{2}\\
&  \leq\mu\otimes\omega\left(  \left[  1_{R}\otimes d\right]  ^{\ast}\left[
1_{R}\otimes d\right]  \right)  \mu\otimes\omega\left(  \alpha(a^{\ast
}a)\right) \\
&  =0
\end{align*}
by (2.2) and since $\omega(a^{\ast}a)=\left\|  \gamma(a)\right\|  ^{2}=0$. But
$\gamma(A)$ is dense in $H$, so $\tilde{\mu}\circ\alpha(a)=0$. Clearly
$\tilde{\mu}^{\alpha}$ is linear, and as in the above calculation we have for
any $a,d\in A$ that%
\[
\left|  \left\langle \gamma(d),\tilde{\mu}^{\alpha}(\gamma(a))\right\rangle
\right|  \leq\mu\left(  1_{R}\right)  \left\|  \gamma(d)\right\|  \left\|
\gamma(a)\right\|
\]
so $\left\|  \tilde{\mu}^{\alpha}\right\|  \leq\mu\left(  1_{R}\right)
=\left\|  \mu\right\|  $. Hence $\tilde{\mu}^{\alpha}$ has a unique bounded
linear extension to $H$, with the same norm. If $\alpha$ is unital, then by
(2.1) we have
\begin{align*}
\left\langle \gamma(d),\tilde{\mu}^{\alpha}\left(  \gamma\left(  1_{A}\right)
\right)  \right\rangle  &  =\left\langle \gamma(d),\tilde{\mu}\left(
1_{R}\otimes1_{A}\right)  \right\rangle \\
&  =\mu\left(  1_{R}\right)  \omega\left(  d^{\ast}1_{A}\right) \\
&  =\left\langle \gamma(d),\mu\left(  1_{R}\right)  \gamma\left(
1_{A}\right)  \right\rangle
\end{align*}
so $\tilde{\mu}^{\alpha}\left(  \gamma\left(  1_{A}\right)  \right)
=\mu\left(  1_{R}\right)  \gamma\left(  1_{A}\right)  $ from which $\left\|
\tilde{\mu}^{\alpha}\right\|  =\left\|  \mu\right\|  $ follows. $\square$

\bigskip

Lastly we will need the following important property in the proof of the mean
ergodic theorem. Note that by a \textit{normal} $\ast$-homomorphism from one
von Neumann algebra to another, we mean a $\ast$-homomorphism that is $\sigma
$-weakly continuous.

\bigskip

\noindent\textbf{Proposition 2.3.}\textit{ Consider the situation in
Propositions 2.1 and 2.2. Furthermore, let }$\nu$\textit{ be another normal
positive linear functional on }$R$\textit{ satisfying }$(\nu\otimes
\omega)\circ\alpha=\nu\left(  1_{R}\right)  \omega$\textit{. Also assume that
}$\alpha$\textit{ is normal, and that }$\Delta:R\rightarrow R\otimes
R$\textit{ is a normal }$\ast$\textit{-homomorphism such that}%
\[
\left(  \iota_{R}\otimes\alpha\right)  \circ\alpha=\left(  \Delta\otimes
\iota_{A}\right)  \circ\alpha
\]
\textit{Write}%
\[
\mu\ast\nu:=(\mu\otimes\nu)\circ\Delta
\]
\textit{then it follows that}%
\[
\widetilde{\mu\ast\nu}^{\alpha}=\tilde{\nu}^{\alpha}\tilde{\mu}^{\alpha}%
\]

\bigskip

\noindent\textbf{Proof.} For any $r\in R$ and $a\in A$ we have for all $d\in
A$ that%
\begin{align*}
\omega\left(  d^{\ast}\left[  \mu\otimes\left(  \nu\otimes\iota_{A}\right)
\right]  \circ\left(  \iota_{R}\otimes\alpha\right)  (r\otimes a)\right)   &
=\omega\left(  d^{\ast}\mu(r)\left(  \nu\otimes\iota_{A}\right)  \circ
\alpha(a)\right) \\
&  =\left\langle \gamma(d),\tilde{\nu}\circ\alpha(\mu(r)a)\right\rangle \\
&  =\left\langle \gamma(d),\tilde{\nu}^{\alpha}(\tilde{\mu}(r\otimes
a))\right\rangle
\end{align*}
by Propositions 2.1 and 2.2, hence by linearity%
\begin{equation}
\omega\left(  d^{\ast}\left[  \mu\otimes\left(  \nu\otimes\iota_{A}\right)
\right]  \circ\left(  \iota_{R}\otimes\alpha\right)  (T)\right)  =\left\langle
\left(  \tilde{\nu}^{\alpha}\right)  ^{\ast}\gamma(d),\tilde{\mu
}(T)\right\rangle \tag{2.3}%
\end{equation}
for all $T\in R\odot A$. The left hand side of (2.3) is a $\sigma$-weakly
continuous linear functional of $T\in R\otimes A$, since $\iota_{R}%
\otimes\alpha$ is the tensor product of two $\sigma$-weakly continuous $\ast
$-homomorphisms, and $\left(  \mu/\left\|  \mu\right\|  \right)
\otimes\left(  \left(  \nu/\left\|  \nu\right\|  \right)  \otimes\iota
_{A}\right)  $ that of two $\sigma$-weakly continuous conditional expectations
(the case $\mu=0$ or $\nu=0$ being trivial). The right hand side of (2.3) is
also a $\sigma$-weakly continuous linear functional of $T\in R\otimes A$. To
see this, consider any net $T_{\lambda}\in R\otimes A$ converging $\sigma
$-weakly to $T$. For any $c\in A$ one has
\[
\left\langle \gamma(c),\tilde{\mu}\left(  T_{\lambda}\right)  \right\rangle
=\mu\otimes\omega\left(  \left[  1_{R}\otimes c\right]  ^{\ast}T_{\lambda
}\right)  \rightarrow\mu\otimes\omega\left(  \left[  1_{R}\otimes c\right]
^{\ast}T\right)  =\left\langle \gamma(c),\tilde{\mu}\left(  T\right)
\right\rangle
\]
in the $\lambda$ limit, since $\mu\otimes\omega$ is $\sigma$-weakly
continuous. However, the $\sigma$-weak topology is a weak* topology, hence by
the resonance theorem (uniform boundedness) the net $\left(  T_{\lambda
}\right)  $ is bounded in the norm of $R\otimes A$. Since $\gamma(A)$ is dense
in $H$, it therefore follows that $\left\langle x,\tilde{\mu}\left(
T_{\lambda}\right)  \right\rangle \rightarrow\left\langle x,\tilde{\mu}\left(
T\right)  \right\rangle $ for all $x\in H$, so indeed (2.3)'s right hand side
is $\sigma$-weakly continuous in $T$. But $R\odot A$ is $\sigma$-weakly dense
in $R\otimes A$, therefore (2.3) holds for all $T\in R\otimes A$, in
particular for $T=\alpha(a)$, so%
\begin{align*}
\left\langle \gamma(d),\tilde{\nu}^{\alpha}\tilde{\mu}^{\alpha}(\gamma
(a))\right\rangle  &  =\omega\left(  d^{\ast}\left[  \mu\otimes\left(
\nu\otimes\iota_{A}\right)  \right]  \circ\left(  \iota_{R}\otimes
\alpha\right)  \circ\alpha(a)\right) \\
&  =\omega\left(  d^{\ast}\left[  (\mu\otimes\nu)\otimes\iota_{A}\right]
\circ\left(  \Delta\otimes\iota_{A}\right)  \circ\alpha(a)\right) \\
&  =\omega\left(  d^{\ast}\left\{  \left[  (\mu\otimes\nu)\circ\Delta\right]
\otimes\iota_{A}\right\}  \circ\alpha(a)\right) \\
&  =\left\langle \gamma(d),\widetilde{\mu\ast\nu}^{\alpha}(\gamma
(a))\right\rangle
\end{align*}
for any $a\in A$, by Propositions 2.1 and 2.2, and since $\Delta$ is normal
(which ensures that $\left[  (\mu\otimes\nu)\otimes\iota_{A}\right]
\circ\left(  \Delta\otimes\iota_{A}\right)  =\left[  (\mu\otimes\nu
)\circ\Delta\right]  \otimes\iota_{A}$ on $R\otimes A$). Since $\gamma(A)$ is
dense in $H$, we obtain $\tilde{\nu}^{\alpha}\tilde{\mu}^{\alpha}%
=\widetilde{\mu\ast\nu}^{\alpha}$. $\square$

\section{The mean ergodic theorem}

Continuing with Section 2's notation, we can now formulate and prove a mean
ergodic theorem:

\bigskip

\noindent\textbf{Theorem 3.1.}\textit{ Consider two normal }$\ast
$\textit{-homomorphisms }$\Delta:R\rightarrow R\otimes R$\textit{ and }%
$\alpha:A\rightarrow R\otimes A$\textit{ such that }$\left(  \iota_{R}%
\otimes\alpha\right)  \circ\alpha=\left(  \Delta\otimes\iota_{A}\right)
\circ\alpha$\textit{ and }$(\theta\otimes\omega)\circ\alpha=\omega$\textit{
for all }$\theta\in\left(  R_{\ast}^{+}\right)  _{1}$\textit{. Assume the
existence of a net }$\left(  \varphi_{\lambda}\right)  $\textit{ in }$\left(
R_{\ast}^{+}\right)  _{1}$\textit{ such that }$\left\|  \theta\ast
\varphi_{\lambda}-\varphi_{\lambda}\right\|  \rightarrow0$\textit{ for all
}$\theta\in\left(  R_{\ast}^{+}\right)  _{1}$\textit{. Let }$P$\textit{ the
projection of }$H$\textit{ on }$V:=\left\{  x\in H:\tilde{\theta}^{\alpha
}x=x\text{\textit{ for all }}\theta\in\left(  R_{\ast}^{+}\right)
_{1}\right\}  $\textit{. Then }%
\[
\lim_{\lambda}\left\langle x,\tilde{\varphi}_{\lambda}^{\alpha}y\right\rangle
=\left\langle x,Py\right\rangle
\]
\textit{for all }$x,y\in H$\textit{.}

\bigskip

\noindent\textbf{Proof.} Set\textit{ }%
\[
N=\overline{\text{span}\left\{  x-\tilde{\theta}^{\alpha}x:x\in H,\theta
\in\left(  R_{\ast}^{+}\right)  _{1}\right\}  }%
\]
and note that $\left\|  \tilde{\theta}^{\alpha}\right\|  \leq\left\|
\theta\right\|  =1$, i.e. $\tilde{\theta}^{\alpha}$ is a contraction, then by
a standard argument $N=V^{\bot}$ (see for example \cite[Section 1.1]{K}). Keep
in mind that $\left(  R\otimes A\right)  _{\ast}=R_{\ast}\otimes_{\ast}%
A_{\ast}$ where by $\otimes_{\ast}$ we mean the tensor product of Banach
spaces with the completion taken in the dual norm of the spatial C*-norm on
$R\odot A$ (see for example \cite[Section 1.22]{Sa}); this will be useful in
the following calculation. Note that this dual norm is a cross norm. For any
$a,d\in A$ and $\theta\in\left(  R_{\ast}^{+}\right)  _{1}$ it follows from
Proposition 2.3 that%
\begin{align*}
&  \left|  \left\langle \gamma(d),\tilde{\varphi}_{\lambda}^{\alpha}\left(
\gamma(a)-\tilde{\theta}^{\alpha}\gamma(a)\right)  \right\rangle \right| \\
&  =\left|  \left\langle \gamma(d),\tilde{\varphi}_{\lambda}^{\alpha}%
\gamma(a)-\widetilde{\theta\ast\varphi_{\lambda}}^{\alpha}\gamma
(a)\right\rangle \right| \\
&  =\left|  \varphi_{\lambda}\otimes\omega\left(  \left[  1_{R}\otimes
d\right]  ^{\ast}\alpha(a)\right)  -\left(  \theta\ast\varphi_{\lambda
}\right)  \otimes\omega\left(  \left[  1_{R}\otimes d\right]  ^{\ast}%
\alpha(a)\right)  \right| \\
&  =\left|  \left(  \varphi_{\lambda}-\theta\ast\varphi_{\lambda}\right)
\otimes\omega\left(  \left[  1_{R}\otimes d\right]  ^{\ast}\alpha(a)\right)
\right| \\
&  \leq\left\|  \varphi_{\lambda}-\theta\ast\varphi_{\lambda}\right\|
\left\|  \omega\right\|  \left\|  \left[  1_{R}\otimes d\right]  ^{\ast}%
\alpha(a)\right\| \\
&  \rightarrow0
\end{align*}
Furthermore $\left\|  \tilde{\varphi}_{\lambda}^{\alpha}-\tilde{\varphi
}_{\lambda}^{\alpha}\tilde{\theta}^{\alpha}\right\|  \leq2$ by Proposition
2.2, so $\left(  \tilde{\varphi}_{\lambda}^{\alpha}-\tilde{\varphi}_{\lambda
}^{\alpha}\tilde{\theta}^{\alpha}\right)  $ is a bounded net, while
$\gamma(A)$ is dense in $H$, hence
\[
\left\langle x,\tilde{\varphi}_{\lambda}^{\alpha}\left(  y-\tilde{\theta
}^{\alpha}y\right)  \right\rangle \rightarrow0
\]
for all $x,y\in H$ and $\theta\in\left(  R_{\ast}^{+}\right)  _{1}$. Since
$\left\|  \tilde{\varphi}_{\lambda}^{\alpha}\right\|  \leq1$, we conclude from
the definition of $N$ that
\[
\left\langle x,\tilde{\varphi}_{\lambda}^{\alpha}y\right\rangle \rightarrow0
\]
for all $x\in H$ and all $y\in N$. So for any $x,y\in H$ we obtain%
\begin{align*}
\left\langle x,\tilde{\varphi}_{\lambda}^{\alpha}y\right\rangle  &
=\left\langle x,\tilde{\varphi}_{\lambda}^{\alpha}Py\right\rangle
+\left\langle x,\tilde{\varphi}_{\lambda}^{\alpha}(1-P)y\right\rangle \\
&  =\left\langle x,Py\right\rangle +\left\langle x,\tilde{\varphi}_{\lambda
}^{\alpha}(1-P)y\right\rangle \\
&  \rightarrow\left\langle x,Py\right\rangle
\end{align*}
by the definition of $P$ and since $(1-P)y\in V^{\bot}=N$. $\square$

\bigskip

In particular this result holds in the situation presented in
Section 1, where $R=M$ is an amenable locally compact quantum group.
This is our main and final result, and we now conclude with a few
brief remarks to give some indication of the relation with classical
ergodic theory.

Note that if $\alpha$ is unital in Theorem 3.1, then one has $P\Omega=\Omega$,
where $\Omega:=\gamma(1_{A})$ is the (non-zero) cyclic vector of $\left(
A,\omega\right)  $'s GNS construction, since
\[
\left\langle \gamma(d),\tilde{\theta}^{\alpha}\Omega\right\rangle
=\left\langle \gamma(d),\tilde{\theta}(1_{R}\otimes1_{A})\right\rangle
=\left\langle \gamma(d),\Omega\right\rangle
\]
for all $\theta\in\left(  R_{\ast}^{+}\right)  _{1}$. This is
essentially the same situation as in classical ergodic theory.
Extending the classical case, it seems reasonable to say that the
\textit{dynamical system} $\left( A,\omega,\alpha\right)  $ is
\textit{ergodic} when $\dim PH=1$, i.e. $PH=\mathbb{C}\Omega$. Using
Theorem 3.1, this is easily seen to be equivalent to
\[
\lim_{\lambda}\varphi_{\lambda}\otimes\omega\left(  \lbrack1_{R}\otimes
a]\alpha(b)\right)  =\omega(a)\omega(b)
\]
again paralleling the situation in classical ergodic theory.

\bigskip

\noindent\textbf{Acknowledgment. }I thank Johan Swart and Gusti van Zyl for
useful conversations.

\end{document}